\documentclass[12pt,twoside]{amsart}
\usepackage{amssymb}

\newtheorem{thm}{Theorem}[section]
\newtheorem{proposition}[thm]{Proposition}

\newtheorem{corollary}[thm]{Corollary}

\theoremstyle{definition}
\newtheorem{definition}[thm]{Definition}
\newtheorem{remark}[thm]{Remark}
\newtheorem{example}[thm]{Example}

\newcommand{\Gin}{\hbox{\rm Gin}}
\newcommand{\reg}{\hbox{\rm reg}}

\begin{document}

\title[The Geometry of Hilbert Functions]{The Geometry of Hilbert
Functions}

\author[Juan C.\ Migliore]{Juan C.\ Migliore}
\address{Department of Mathematics,
            University of Notre Dame,
            Notre Dame, IN 46556,
            USA}
\email{migliore.1@nd.edu}

\maketitle

\section{Introduction} 

The title of this paper, ``The geometry of Hilbert functions," might
better be suited for a multi-volume treatise than for a
single short article.  Indeed, a large part of the beauty of, and interest in,
Hilbert functions derives from their ubiquity in all of commutative
algebra and algebraic geometry, and the unexpected information that they
can give, very much of it expressible in a geometric way.  Most of this paper is
devoted to describing just one small facet of this theory, which connects results of
Davis (e.g.\ \cite{davis}) in the 1980's, of Bigatti, Geramita and myself
(cf.\ \cite{BGM}) in the 1990's, and of Ahn and myself
(cf.\ \cite{AM}) very recently.  On the other hand, we have an alphabet soup of
topics that play a role here: UPP, WLP, SLP, ACM at the very least.  It is
interesting to see the ways in which these properties interact, and we also try to
illustrate some aspects of this.

There are almost as many different notations for Hilbert functions as
there are papers on the subject.  We will use the following notation.  If
$I$ is a homogeneous ideal in a polynomial ring $R$, we write
\[
h_{R/I}(t) := \dim (R/I)_t.
\]
If $I$ is a saturated ideal defining a subscheme $V$ of $\mathbb P^n$
then we also write this function as $h_V(t)$ or $h_{R/I_V}(t)$.  

So where is the geometry?  Of course $\dim R_t = \binom{t+n}{n}$, so the information
provided by the Hilbert function is equivalent to giving the dimension of the degree
$t$ component of $I$.  This dimension is one more than the dimension of the linear
system of hypersurfaces of degree $t$ defined by $I_t$ (since this latter dimension
is projective).  What is the base locus of this linear system?  Of course $V$ is
contained in this base locus, but it may contain more.  The results in this paper
(e.g.\ Theorem \ref{BGM general results}, Theorem \ref{AM general}, Theorem \ref{BGM
UPP results} and Theorem \ref{AM UPP results}) can be viewed as describing the
dimension, irreducibility and reducedness of this base locus, based on information
about the Hilbert function, and other basic properties, of
$V$.  We will see that under some situations, just knowing the dimension of this
linear system in two consecutive degrees can force the base locus to contain a
hypersurface, or anything smaller.  (We will concentrate on the
 curve case.)

An important starting point for us (and indeed for almost any discussion of Hilbert
functions of standard graded algebras) is Macaulay's theorem bounding the growth of
the Hilbert function.  Once we have this, we need Gotzmann's results about what
happens when Macaulay's bound is achieved.  These are both discussed in Section
\ref{preliminary section}, as are several other results related to these.

In Section \ref{WLP section} we recall the notions of the Uniform Position Property
(UPP) and the Weak Lefschetz Property (WLP) and some of their connections. 
Subsequent sections, especially Section \ref{UPP results}, continue the discussion
of UPP.  WLP, while often less visible, lurks in the background of many
of the results and computations of this paper, and in fact is an important object of
study.  We include a short discussion of the behavior of WLP in families of points
in Section \ref{WLP section}, including a new example (Example \ref{WLP in
families}) showing how, for fixed Hilbert function, WLP can hold in one component of
the postulation Hilbert scheme and not hold in another.  See also Theorem~\ref{delta
2}.

The focus in this article is the situation where the first difference of the Hilbert
function of a set of points, $Z$, in projective space $\mathbb P^n$ attains the same
value in two consecutive degrees: $\Delta h_Z(d) = \Delta h_Z(d+1) = s$.  Depending
on the relation between $d$, $s$ and certain invariants of $Z$, we will get
geometric consequences for the base locus.  In Section \ref{set stage} we describe
these relations, setting the stage for the main results. 

These main results are given in Sections \ref{general results} and \ref{UPP
results}.  Here we see that under certain assumptions on $d$, the condition $\Delta
h_Z(d) = \Delta h_Z(d+1) = s$ guarantees that the base locus of the linear system
$|I_d|$ is a curve of degree $s$.  This comes from work  in \cite{davis}, \cite{BGM}
and \cite{AM}.  Other results follow as well.  What is surprising here is that the
central condition of \cite{AM}, namely that $d > r_2(R/I_Z)$ (see Section
\ref{preliminary section} for the definition), is much weaker than the central
assumption of the comparable results in  [Bigatti-Geramita-Migliore],
namely
$d
\geq s$, but the results are very similar.  Section \ref{general results} focuses on
the general results, while Section \ref{UPP results} turns to the question of what
can be said about this base locus when the points have UPP.

There are some differences in the results of [Bigatti-Geramita- \linebreak Migliore]
and
\cite{AM} as a result of the differences in these assumptions.  Section
\ref{example section} studies these, and gives examples to show that they are not
accidental omissions.  Some very surprising behavior is exhibited here.

  I am grateful to Irena Peeva for asking me to write this paper, which I
enjoyed doing.  In part it is a greatly expanded version of a talk that I
gave in the Algebraic Geometry seminar at Queen's University in the fall of 2004,
and I am grateful to Mike Roth and to Greg Smith for their kind invitation.  I would
like to thank Jeaman Ahn, Chris Francisco, Hal Schenck and especially Tony Iarrobino
for helpful comments.  And of course I am most grateful to my co-authors Anna Bigatti
and Tony Geramita (\cite{BGM}) and Jeaman Ahn (\cite{AM}) for their insights and for
the enjoyable times that we spent in our collaboration. During the writing of this
paper, and some of the work described here, I was sponsored by the National Security
Agency (USA)  under Grant Number MDA904-03-1-0071.


\section{Maximal growth of the Hilbert function} \label{preliminary section}

We first collect the notation that we will use throughout this paper.  Let $k$ be a
field of characteristic zero and let
$R = k[x_1,\dots, x_n]$.  

\begin{definition} 
Let $Z \subset \mathbb P^{n-1}$ be any closed subscheme with defining (saturated)
ideal
$I = I_Z$.  

\begin{itemize}
\item[(a)] The {\em Hilbert function of $Z$} is the function 
\[
h_Z(t) = \dim(R/I_Z)_t
\]
We also may write $h_{R/I}(t)$ for this function. 
If $A$ is Artinian then we write
\[
h_A(t) = \dim A_t
\]
for its Hilbert function.

\item[(b)] We say that $Z$ is {\em arithmetically Cohen-Macaulay} (ACM) if the
coordinate ring $R/I_Z$ is a Cohen-Macaulay ring.  Note that if $Z$ is a
zero-dimensional scheme then it is automatically ACM. \qed
\end{itemize}
\end{definition}

If $F$ is a homogeneous polynomial, by abuse of notation we will also denote by $F$
the hypersurface of $\mathbb P^{n-1}$ defined by $F$.

\begin{definition} \label{alpha definition}
 For a homogeneous ideal $I$ we define
\[
\alpha = \min \{ t \ | \ I_t \neq 0 \},
\]
i.e.\ $\alpha$ is the {\em initial degree} of $I$. \qed
\end{definition}

If $A = R/I$ is a standard graded $k$-algebra, then there is a famous bound, due to
Macaulay (cf.\ \cite{fsmacaulay}), that describes the maximum possible growth of the
Hilbert function of $A$ from any degree to the next.  To give this bound, we need a
little preparation.  

\begin{definition} \label{ibinomexp}
The {\em $i$-binomial expansion} of the integer $c$ ($i, c >0$) is the unique
expression 
\[
c = \binom{m_i}{i} + \binom{m_{i-1}}{i-1} + \dots + \binom{m_j}{j},
\]
where $m_i > m_{i-1} > \dots > m_j \geq j \geq 1$.  \qed
\end{definition}

Note that the assertion that this representation is unique is something that has to
be checked!

\begin{definition}
If $c \in {\mathbb Z}$ ($c>0$) has $i$-binomial expansion as in Definition
\ref{ibinomexp}, then we set  
\[
c^{\langle i \rangle} = \binom{m_i+1}{i+1} + \binom{m_{i-1}+1}{i} + \dots +
\binom{m_j+1}{j+1}.
\]
Note that this defines a collection of functions $^{\langle i \rangle}: \mathbb Z
\rightarrow \mathbb Z$.  \qed
\end{definition}

For example, the 5-binomial expansion of 76 is
\[
76 = \binom{8}{5}  + \binom{6}{4} + \binom{4}{3} + \binom{2}{2},
\]
so
\[
76^{\langle 5 \rangle} = \binom{9}{6} + \binom{7}{5} + \binom{5}{4} + \binom{3}{3} =
111.
\]

\begin{definition} \label{o-seq}
A sequence of non-negative integers $\{ c_i : i \geq 0 \}$ is called an {\em
O-sequence} if
\[
c_0 = 1 \ \hbox{ and } \ c_{i+1} \leq c_i^{\langle i \rangle}, 
\]
for all $i$.  An O-sequence is said to have {\em maximal growth from degree $i$ to
degree $i+1$} if $c_{i+1} = c_i^{\langle i \rangle}$.  \qed
\end{definition}

The importance of the binomial expansions described above becomes apparent from the
following beautiful theorem of Macaulay:

\begin{thm}[\cite{fsmacaulay}] \label{macaulay thm}

The following are equivalent:

\begin{itemize}
\item[(i)] $\{ c_i : i \geq 0 \}$ is an O-sequence;

\item[(ii)] $\{ c_i : i \geq 0 \}$ is the Hilbert function of a standard graded
$k$-algebra.
\end{itemize}
\end{thm}

\noindent In other words, a sequence of non-negative integers is the Hilbert
function of a standard graded $k$-algebra if and only if the growth from any degree
to the next is bounded as above.  Remember that the Hilbert function is eventually
equal to the Hilbert polynomial, at which point $c_{i+1} = c_i^{\langle i \rangle}$,
but until that point, anything is allowed as long as we maintain $c_{i+1} \leq 
c_i^{\langle i \rangle}$.

One can ask when such a sequence is the Hilbert function of a {\em reduced}
$k$-algebra.  This was answered by Geramita, Maroscia and Roberts (cf.\ \cite{GMR}),
by introducing the {\em first difference} of the Hilbert function.

\begin{definition}
Given a sequence of non-negative integers $ \underline{c} = $\linebreak $\{ c_i : i
\geq 0
\}$, the {\em first difference} of this sequence is the sequence $\Delta
\underline{c} :=
\{ b_i \}$ defined by $b_i = c_i - c_{i-1}$ for all $i$.  (We make the convention
that $c_{-1} = 0$, so $b_0 = c_0 = 1$.)  We say that $\underline{c}$ is a {\em
differentiable O-sequence} if $\Delta \underline{c}$ is again an O-sequence.  By
taking successive first differences, we inductively define the {\em $k$-th
difference}
$\Delta^k \underline{c}$.
\qed
\end{definition}

\begin{remark} \label{first diff} An important fact to remember is that if $Z$ is a
zero-dimensional scheme with Artinian reduction $A$ then
\[
h_A(t) = \Delta h_Z(t) \hbox{ for all $t$}.
\]
It follows that  ($\Delta h_Z(t) : t \geq 0)$ is a
finite sequence of positive integers, called the {\em $h$-vector} of $Z$. 
Similarly, if $V$ is ACM of dimension $d$ with Artinian reduction $A$ then
$h_A(t) = \Delta^{d+1} h_V (t) \hbox{ for all $t$}$.  \qed
\end{remark}

The following is the theorem of Geramita, Maroscia and Roberts mentioned above.  It
guarantees the existence of a reduced subscheme of $\mathbb P^{n-1}$ with given
Hilbert function (remember that
$R = k[x_1,\dots,x_n]$), under a simple hypothesis on the Hilbert function.

\begin{thm}[\cite{GMR}] \label{GMR theorem}
Let $ \underline{c} = \{ c_i \}$ be a sequence of non-negative integers, with $c_1 =
n$.  Then $\underline{c}$ is the Hilbert function of a standard graded $k$-algebra
$R/I$, with $I$ radical, if and only if $\underline{c}$ is a differentiable
O-sequence. 
\end{thm}

Note that if $I$ is any saturated ideal (reduced or otherwise), and if $L$ is a
general linear form, then we have $[I:L] \cong I(-1)$ as graded modules.  It follows
that $L$ induces an injection $ \times L : ((R/I)(-1))_t \rightarrow (R/I)_t$, i.e.\
$R/I$ has depth at least 1.  Hence the first difference of the Hilbert function of
$R/I$ is the Hilbert function of $R/(I,L)$, again a standard graded $k$-algebra,
so the Hilbert function of $R/I$ is a differentiable O-sequence.  This  shows that
Theorem \ref{GMR theorem} would also be true if ``radical" were replaced by
``saturated" (since if a radical ideal can be constructed for a given
differentiable O-sequence, this ideal is also saturated).  The real heart of Theorem
\ref{GMR theorem} is that any such sequence can be achieved by a radical ideal.

\begin{example}
We should stress that Theorem \ref{GMR theorem} guarantees the existence of a reduced
(even non-degenerate) subscheme of $\mathbb P^n$, but it does not guarantee that
what we get will be equidimensional.  It also does not say anything about higher
differences.  (Note that even if $I$ is saturated, this is not necessarily true of
$(I,L)$.)  We give two examples.

\begin{enumerate}
\item Let $I = (x_1,x_2)$ in $k[x_0,x_1,x_2,x_3]$, so $I$ defines a line, $\lambda$,
in $\mathbb P^3$.  Let $F$ be a homogeneous polynomial of degree 3, non-singular
along $\lambda$.  Let $J = (I^2,F)$.  Then $J$ is the saturated ideal of a
non-reduced subscheme of degree 2 and genus $-2$ (cf.\ \cite{dble}).  Its Hilbert
function is the sequence $1, 4, 7, 9, 11, \dots$ (with Hilbert polynomial $2t+3$). 
The smallest genus for a reduced equidimensional subscheme of degree 2 is $-1$ (two
skew lines), so this Hilbert function does not exist among reduced, equidimensional
curves in $\mathbb P^3$.  However, adding points reduces the genus while not
affecting the degree, and in fact this Hilbert function occurs for the union of two
skew lines and one point, which is indeed reduced.

\item Consider the sequence $\underline{c} = 1, 4, 10, 17, 26, 35, \dots $ (with
Hilbert polynomial $9t - 11 +1$, so it corresponds to a curve of degree 9 and
arithmetic genus 11).  Its first difference is $1, 3, 6, 7, 9, 9, \dots$, which is
again an O-sequence, so $\underline{c}$ is a differentiable O-sequence.  However,
the second difference is $1, 2, 3, 1, 2$ which is not an O-sequence.  Theorem
\ref{GMR theorem} guarantees the existence of a reduced curve with Hilbert function
$\underline{c}$, and indeed it can be achieved by the union $C_1 \cup C_2 \cup P_1
\cup P_2$, where $C_1$ is a plane curve of degree~3, $C_2$ is a plane curve of
degree 6 not in the plane of $C_1$, $C_1$ and $C_2$ meet in 3 points, $P_1$ is a
generally chosen point in $\mathbb P^3$, and $P_2$ is a generally chosen point in the
plane of $C_1$. 
\end{enumerate}
 One sees that finding the reduced subscheme of Theorem \ref{GMR
theorem} can be a tricky matter!  \qed
\end{example}

If $\underline{h} = \{ c_i : i \geq 0 \}$ is an O-sequence, Theorem \ref{macaulay
thm} guarantees that $\underline{h}$ is the Hilbert function of a standard graded
$k$-algebra, say $R/I$.  The following striking result of Gotzmann says what happens
if we have maximal growth (see Definition \ref{o-seq}), i.e.\ if the equality
$c_{d+1} = c_d^{\langle d \rangle}$ holds.  Let $\bar I$ be the ideal generated by
the components of $I$ of degree $\leq d$, so that $\bar I_i = I_i$ for all $i \leq
d$.  We also sometimes write $\bar I = \langle I_{\leq d} \rangle$.

\begin{thm}[\cite{gotzmann}, \cite{green}] \label{gotzmann persistence}
Let $\underline{h} = \{ c_i : i \geq 0 \}$ be the Hilbert function of $R/\bar I$, and
suppose that $c_d = \binom{m_d}{d} + \binom{m_{d-1}}{d-1} + \dots + \binom{m_j}{j}$
is the $d$-binomial expansion of $c_d$. Assume that $c_{d+1} = c_d^{\langle d
\rangle}$.  Then for any
$l
\geq 1$ we have
\[
c_{d+l} = \binom{m_d +l}{d+l} + \binom{m_{d-1}+l}{d-1+l} + \dots + \binom{m_j
+l}{j+l}.
\]
In particular, the Hilbert polynomial $p_{R/\bar I}(x)$ agrees with the Hilbert
function $h_{R/\bar I}(x)$ in all degrees $\geq d$.  It can be written as
\[
p_{R/\bar I}(x) = \binom{x+m_d-d}{m_d-d} + \binom{x+m_{d-1}-d}{m_{d-1}-(d-1)} + \dots
+
\binom{x+m_j-d}{m_j-j},
\]
where $m_d-d \geq m_{d-1}-(d-1) \leq \dots \leq m_j-j$.
\end{thm}

\begin{remark} \label{interpret max gr} Theorem \ref{gotzmann persistence} is also
known as the {\em Gotzmann Persistence Theorem}.  Since $p_{R/\bar I}(x)$ has degree
$m_d-d$, this means that the Krull dimension $\dim R/\bar I = m_d-d+1$, and $\bar I$
defines a subscheme of dimension
$m_d-d$.  Its degree can also easily be computed by checking the leading
coefficient.  Gotzmann also shows that $\bar I$ has regularity $\leq d$ and agrees
with its saturation in degrees $\geq d$.
\qed
\end{remark}

\begin{remark} \label{gotzmann references}
An excellent reference for the Gotzmann Persistence Theorem and related results is
section 4.3 of  \cite{bruns-herzog2} (note that this is the revised edition; this
material does not appear in the original book).  A {\em much} more detailed and
comprehensive expository treatment  than that given here can be found in
\cite{iarrobino-kanev}, Appendix C (written by A.\ Iarrobino and S.\ Kleiman).  This
includes a detailed description of the associated Hilbert scheme.  \qed
\end{remark}

Much of the work described in this paper revolves around the following idea.  If $I$
is a saturated homogeneous ideal and $L$ is a general linear form, then reducing
modulo $L$ gives a new ideal, $J$, in $S=R/(L)$.  Suppose that $S/J$ has maximal
growth from degree $d$ to degree $d+1$.  This does not necessarily imply that $R/I$
has maximal growth at that spot (see also Remark \ref{delta max not imply I max}). 
Then what information {\em can} we get about $I$ if we know that $R/J$ has maximal
growth?  As a first step we have the following result from \cite{BGM}.  We have
noted in Remark \ref{interpret max gr} that maximal growth of $J$ implies something
about the saturation and regularity of $J$, thanks to Gotzmann, but we now show that
it also says something about the saturation and regularity of $\langle I_{\leq d}
\rangle$.

\begin{proposition}{\rm ([Bigatti-Geramita-Migliore] Lemma 1.4, Proposition 1.6)}
\label{BGM saturation result} Let $I \subset R$ be a saturated ideal.  Let $L$ be a
general linear form, and let $J = \frac{(I,L)}{(L)} \subset S = R/(L)$.  Suppose
that $S/J$ has maximal growth from degree $d$ to degree $d+1$.  Let $\bar I =
\langle I_{\leq d} \rangle$ as above.  Let $\bar I^{sat}$ be the saturation of $\bar
I$.  Then $\bar I =
\bar I^{sat}$, i.e.\ $\bar I$ is a saturated ideal.  Furthermore, $\bar I$ is
$d$-regular.
\end{proposition}

\begin{remark} \label{comparison}
Again, Proposition \ref{BGM saturation result} differs from the Gotzmann results
(cf.\ Theorem \ref{gotzmann persistence}, Remark \ref{interpret max gr}) in that
$S/J$ is assumed to have maximal growth, but we conclude something about
$I$.  But even beyond this, there is another difference between Proposition \ref{BGM
saturation result} and the part of Gotzmann's work that we have presented (referring
to the sources mentioned in Remark \ref{gotzmann references} for a more complete
exposition); indeed, in Theorem \ref{gotzmann persistence} and Remark \ref{interpret
max gr} we only conclude that the ideal in question agrees with its saturation in
degree $d$ and beyond, while here we make a conclusion about the whole ideal. (The
price is that we have to assume that $I$ itself is saturated to begin with,
although this does not necessarily hold for $J$.)

Also, the maximal growth assumption for $S/J$ is necessary.  For instance, let $I$
be the saturated ideal of a set, $Z$, of sixteen general points in $\mathbb P^3$. 
Then $I$ has four generators in degree 3 and three generators in degree 4.  The
four generators in degree 3 do define $Z$ scheme-theoretically, but this is not
enough.  One can check (e.g.\ using {\tt macaulay} \cite{macaulay}) that if $\bar I =
\langle I_{\leq 3} \rangle$, then  $\bar I^{sat} = I \neq \langle I_{\leq 3} \rangle
= \bar I$.  Furthermore, $\bar I$ is 5-regular and $I$ itself is 4-regular, and
neither is 3-regular.  And indeed, the Hilbert function of $S/J$ is $1,3,6,6,0$, so
we do not have maximal growth of $S/J$ from degree 3 to degree 4. 

Despite this, the results in \cite{AM} show that statements along the lines of
Proposition \ref{BGM saturation result} can be obtained even weakening the maximal
growth assumption.
\qed
\end{remark}


\section{UPP and WLP}  \label{WLP section}

A very important property of (many) reduced sets of points in projective space is
the following:

\begin{definition}  \label{UPP definition}
 A reduced set of points $Z \subset \mathbb P^{n-1}$ has the {\em
Uniform Position Property (UPP)} if, for any $t \leq |Z|$, all subsets of $t$ points
have the same Hilbert function, which necessarily is the truncated Hilbert
function.  \qed
\end{definition}

The last comment follows from the fact that it was also shown in \cite{GMR} that
given any reduced set $Z$ of, say, $d$ points with known Hilbert function, and
truncating this function at any value $k <d$, there is a subset $X$ of $Z$
consisting of
$k$ points whose Hilbert function is this truncated function.  Hence if $Z$ has UPP,
then {\em all} subsets have this truncated Hilbert function.

While it is known that the general hyperplane section of an irreducible curve has
UPP (cf.\ \cite{HE}), at least in characteristic zero, much remains open.  An
important open question is to determine all possible Hilbert functions of sets of
points in projective space with UPP.  It is known in $\mathbb P^2$ but it is open
even in
$\mathbb P^3$.  See
\cite{GM3} for a discussion of several of the many papers that have contributed to
this question for $\mathbb P^2$.  

For any given Hilbert function, there may or may not be a set of points with UPP
having that Hilbert function.  One of the contributions of the papers discussed here
is toward showing some conditions on the Hilbert function that prohibit the
existence of points with UPP having that function.  If there {\em is} a set of
points with UPP having a given Hilbert function, somehow ``most" sets of points with
that Hilbert function have UPP: consider the postulation Hilbert scheme
parameterizing  sets of points with that Hilbert function; then in the component 
containing the given set of points, there is an open subset corresponding to points
with UPP.  (We make this more precise in the discussion about WLP below.)

Another important property is the following. 

\begin{definition} \label{WLP definition}
 An Artinian algebra $A$ has the {\em Weak Lefschetz Property (WLP)} if, for a
general linear form $L$, the map $\times L : A_t \rightarrow A_{t+1}$ has maximal
rank, for all $t$.  We say that $A$ has the {\em Strong Lefschetz Property (SLP)}
if for every $d$, and for a general form $F$ of degree $d$, the map $\times F : A_t
\rightarrow A_{t+d}$ has maximal rank, for all $t$.  If $Z$ is a set of points, we
sometimes say that $Z$ has WLP or SLP if its general Artinian reduction does
\qed
\end{definition}

\begin{remark} 
We have the following comments about WLP and SLP:

\begin{enumerate}

\item The statement that SLP holds for an ideal of any number of general
forms is equivalent to the well-known Fr\"oberg conjecture.  See
\cite{anick}, \cite{MMR2},
\cite{MMR3}, \cite{MMRN}.

\item  {\em Every} Artinian
complete intersection in $k[x_1,x_2,x_3]$ has WLP
(cf.\ \cite{HMNW}).  It is also known that SLP (and hence WLP) holds for {\em every}
ideal in  $k[x_1,x_2]$ (cf.\ \cite{iarrobino}, \cite{HMNW}), and it  (and hence
Fr\"oberg's conjecture) holds for an ideal of general forms in 
$k[x_1,x_2,x_3]$ (cf.\ \cite{anick}).  \qed

\end{enumerate}
\end{remark}

We now begin a short digression about the behavior of WLP in families of
reduced zero-dimensional schemes.
Fix a Hilbert function, $H$, that corresponds to a zero-dimensional
scheme, and consider the $h$-vector, $\underline{h} = 
(a_0,a_1,a_2,\dots,a_s)$, associated to that Hilbert function (i.e.\ its first
difference $\Delta H$ -- cf.\ Remark \ref{first diff}).  Let $d = \sum_i a_i$ be the
degree of the zero-dimensional scheme.  Consider the postulation Hilbert
scheme, ${\mathcal H}_{\underline{h}} := \hbox{Hilb}^{H}(\mathbb P^n)$,
parameterizing  zero-dimensional schemes in $\mathbb P^n$ with  Hilbert function
$H$ (inside the punctual Hilbert scheme, $\hbox{Hilb}^d(\mathbb P^n)$,
of all zero-dimensional schemes in $\mathbb P^n$ with the given degree).  It 
is known that the closure, $\overline{ {\mathcal H}}_{\underline{h}}$,
of ${\mathcal H}_{\underline{h}}$ may have several irreducible
components  -- see for instance Richert (cf.\ \cite{richert}, and use
Hartshorne's lifting procedure, cf.\ \cite{hartshorne}, \cite{MN2}, on the Artinian
monomial ideals that Richert gives), Ragusa-Zappal\`a (cf.\ \cite{RZ}) or Kleppe
(cf.\ \cite{kleppe}, Remark 27).  We will consider only those components of 
$\overline{ {\mathcal H}}_{\underline{h}}$ for which the general element is reduced.

One can show that, like UPP, WLP is an open condition in the sense that in any
component of $\overline{ {\mathcal H}}_{\underline{h}}$, an open
subset (possibly empty) corresponds to zero-dimensional schemes with WLP.  We will
now give an example that answers (positively) the following question.  Namely, does
there exist a Hilbert function with $h$-vector
$\underline{h}$, and components
${\mathcal H}_1$ and ${\mathcal H}_2$ of the corresponding $\overline{{\mathcal
H}}_{\underline{h}}$, for which

\begin{itemize}
\item $\underline{h} =  (a_0,a_1,a_2,\dots,a_s)$ is unimodal, and in fact
satisfies
\begin{equation} \label{unimodal plus}
a_0 < a_1 < a_2 \dots < a_t \geq a_{t+1} \geq \dots \geq a_s.
\end{equation}
for some $t$ (this is a technical necessity for WLP-- cf.\
[Harima-Migliore-Nagel-Watanabe], Remark~3.3);

\item the general element of ${\mathcal H}_1$ corresponds to a reduced,
zero-dimen\-sional scheme with WLP; 

\item {\em no} element of ${\mathcal H}_2$ corresponds to a zero-dimensional scheme
with WLP (i.e.\ the open subset referred to above is empty)?
\end{itemize}

\medskip

The following example answers this question.

\begin{example} \label{WLP in families}

We will give an example of ideals $I_1$ and $I_2$ of reduced zero-dimensional
schemes in $\mathbb P^3$, both with $h$-vector 
\[
\underline{h} = (1,3,6,9,11,11,11),
\]
such that  the Artinian reduction of $R/I_1$ has WLP and the Artinian reduction of
$R/I_2$ does not.  Furthermore, the Betti diagrams for $R/I_1$ and $R/I_2$ are
(respectively)

{\footnotesize 
\begin{verbatim}
total:      1    14    24    11      total:      1    14    25    12 
--------------------------------     --------------------------------
    0:      1     -     -     -          0:      1     -     -     - 
    1:      -     -     -     -          1:      -     -     -     - 
    2:      -     1     -     -          2:      -     1     -     - 
    3:      -     1     -     -          3:      -     1     1     - 
    4:      -     1     2     -          4:      -     2     2     1 
    5:      -     -     -     -          5:      -     -     -     - 
    6:      -    11    22    11          6:      -    10    22    11 
\end{verbatim}
}

Clearly these Betti diagrams allow no minimal element in the sense of Richert
(cf.\ \cite{richert}) or Ragusa-Zappal\`a (cf.\ \cite{RZ}).  Furthermore, neither can
be a specialization of the other, so they correspond to different components of 
$\overline{{\mathcal H}}_{\underline{h}}$.

For $I_1$ we begin with a line in $\mathbb P^3$ and link, using a complete
intersection of type $(3,4)$, to a smooth curve, $C$, of degree 11.  Since a line is
ACM, the same is true of $C$ by liaison.  We let
$Z_1$ be a set of 52 general points on $C$, and let $I_1$ be the homogeneous ideal of
$Z_1$.   It is easy to check that $R/I_1$ has the desired $h$-vector, using liaison
computations (see for instance \cite{migbook}) to compute the Hilbert function of
$C$ and then the fact that $Z_1$ is chosen generically on $C$ so that the Hilbert
function is the truncation.  The rows 0 to 5 in the Betti diagram come only from $C$,
and then the last row is forced from the Hilbert function.  Because the Hilbert
function of $Z_1$ agrees with that of $C$ up to and including degree 6, the WLP for
$Z_1$ follows from the Cohen-Macaulayness of $C$.  Hence the general element of the
corresponding component has WLP.

For $I_2$ we start with the ring $R = k[x,y,z]$ (dropping
subscripts on the variables for convenience).  Consider the monomial ideal $J$
consisting of $(x^3, x^2y^2, x^2yz^2, z^5)$ together with all the monomials of degree
7.  One can check on {\tt macaulay} (cf.\ \cite{macaulay}) that $R/J$ is Artinian
with Hilbert function
$\underline{h}$ and with the Betti diagram above to the right.  Using the lifting
procedure for monomial ideals (cf.\ \cite{hartshorne}, \cite{MN2}), we lift $J$ to
the ideal $I_2$ of a reduced zero-dimensional scheme, $Z_2$.  Now, it is clear from
the Betti diagram that $R/J$ has a socle element in degree 4.  This means that the
map from $(R/J)_4$ to $(R/J)_5$ (both of which are 11-dimensional) induced by a
general linear form necessarily has a kernel, and so it is neither injective nor
surjective.  Hence $Z_2$ does not have WLP.  But from the Betti diagram, it is clear
from semicontinuity that the general element (hence every element) of the component
corresponding to $Z_2$ similarly fails to have WLP.

Incidentally, the lex-segment ideal corresponding to this Hilbert function (and hence
having maximal Betti numbers) has Betti diagram

\begin{verbatim}
             total:      1    18    31    14 
             --------------------------------
                 0:      1     -     -     - 
                 1:      -     -     -     - 
                 2:      -     1     -     - 
                 3:      -     1     1     - 
                 4:      -     2     4     2 
                 5:      -     2     3     1 
                 6:      -    12    23    11 
\end{verbatim}

One can check that this is a specialization of both of the Betti diagrams above.

It would be nice to find an example where $\underline{h}$ is actually the $h$-vector
of a complete intersection.
We remark that it is possible to have a reduced zero-dimensional scheme with the
Hilbert function of a complete intersection, that does not have WLP.  For example,
for the $h$-vector $(1,3,4,4,3,1)$, which is that of a complete intersection of type
$(2,2,4)$, one can take the union, $Z$, of a general set of points in the plane with
$h$-vector $(1,2,3,4,3,1)$ (easily produced with liaison) and two general points in
$\mathbb P^3$, to produce a set of points with the desired $h$-vector.  One can check
that the Artinian reduction of $R/I_Z$ has a socle element in degree 2, so the
multiplication from the component in degree 2 to the component in degree 3 (both of
which are 4-dimensional) induced by a general linear form has no chance to be
injective or surjective, i.e.\ $Z$ does not have WLP.  We know that a complete
intersection in $\mathbb P^3$ has Artinian reduction with WLP (cf.\ \cite{HMNW}), but
we do not know if this example is contained in the same component as the complete
intersection or not.  \qed
\end{example}

\begin{remark}
In the preceding example, note that $Z_1$ has UPP, while $Z_2$ is very far from
having UPP, being the lifting of a monomial ideal.  This motivates the following
natural questions:

\begin{enumerate}
\item Does a set of points with UPP automatically have WLP? 

\item Does the general hypersurface section of a smooth curve necessarily have WLP
(we know that it does have UPP, at least in characteristic zero)?

\item Does the general hyperplane section of a smooth curve necessarily have WLP?  Of
course if the curve is in $\mathbb P^3$ then the hyperplane section is in $\mathbb
P^2$, so the Artinian reduction is an ideal in $k[x_1,x_2]$, and WLP and SLP are
automatic (as we mentioned above).  So the first place that this question is
interesting is for a smooth curve in $\mathbb P^4$.

\item[4.] Does the Artinian reduction of every reduced, arithmetically \linebreak
Gorenstein set of points  have WLP?

\begin{itemize}
\item It is not true that every Artinian Gorenstein ideal in codimension $\geq 5$
has WLP, because the $h$-vector can fail to satisfy condition (\ref{unimodal plus})
on page \pageref{unimodal plus} (cf.\ for instance \cite{ber-iar}, \cite{BL}).  But
the question is open for the Artinian reduction of arithmetically Gorenstein points
in any codimension.  (Part of what is missing is knowledge of what Artinian algebras
lift to reduced sets of points.)

\item It is an open question whether all height 3 Artinian Gorenstein ideals have
WLP.  It is true for all height 3 complete intersections (cf.\ \cite{HMNW}). 
\end{itemize}

\end{enumerate}

In \cite{AM} we were able to answer the first three of the above questions in the
negative, essentially with one example (cf.\ [Ahn- \linebreak Migliore], Examples 6.9
and 6.10).  We omit the details here, but the basic idea is as follows.  We consider
a smooth arithmetically Buchsbaum curve whose deficiency module 
\[
M(C) = \bigoplus_{t \in {\mathbb Z}} H^1(\mathbb P^3, {\mathcal I}_C (t))
\]
 is one-dimensional in degrees 3 and 4. Any linear form induces a homomorphism from
$H^1(\mathbb P^3, {\mathcal I}_C(3))$ to $H^1(\mathbb P^3, {\mathcal I}_C(4))$.  The
condition ``arithmetically Buchsbaum" means that this map is zero for all linear
forms $L$.  It is known that such a smooth curve exists (cf.\ \cite{BM2}).  Now, we
put a lot of points on such a curve, placed generically, and call the resulting set
of points $Z$.  Note that $I_Z$ agrees with $I_C$ in low degrees.  We then 
play with the coordinate ring of $C$, the coordinate ring of $Z$, the
cohomology of ${\mathcal I}_C$ and the cohomology of ${\mathcal I}_Z$, and in
the end we show that the Artinian reduction of $R/I_Z$ does not have WLP. 
Tweaking this example somewhat (looking at the cone over $C$ and taking a general
hypersurface and hyperplane section) answers the remaining questions.  \qed
\end{remark}

Some of the results given below, in particular those coming from the paper
\cite{AM},  are given in terms of {\em reduction numbers}, which we now define.

\begin{definition}[adaptation of Hoa-Trung \cite{hoa-trung}]

Let $I \subset R$ be a homogeneous ideal.  Let $m \geq \dim R/I$.  The $m$-reduction
number of $R/I$ is
\[
\begin{array}{rcl}
r_m(R/I)  & = & \min \{ k \ | \ x_{n-m}^{k+1} \in \Gin(I) \} \\
& = & \min \{ k \ | \ h_{R/(I+J)} \hbox{ vanishes in degree } k+1 \}
\end{array}
\]
where $J$ is an ideal generated by $m$ general linear forms.  We will use the second
line as the definition, and we include the first line only for completeness.  (For
the definition of, and results on, the generic initial ideal, $\Gin(I)$, see for
instance
\cite{green2}.)
\qed
\end{definition}

\begin{example} \label{334 example}
Let $Z$ be a complete intersection of type $(3,3,4)$ in $\mathbb P^{3}$ and let $I =
I_Z$.  Note $\dim R/I = 1$.  Let $L_1, L_2, L_3$ be general linear forms. It is known
that any Artinian reduction of $R/I$ has WLP (cf.\ \cite{HMNW}).

 We have
\[
\begin{array}{rclcll}
h_{R/(I+(L_1))} & : & 1 \ 3 \ 6 \ 8 \ 8 \ 6 \ 3 \ 1 & \Rightarrow & r_1(R/I) = 7. \\

h_{R/(I+(L_1,L_2))} & : & 1 \ 2 \ 3 \ 2 & \Rightarrow & r_2(R/I) = 3 & \hbox{(by
WLP)} \\

h_{R/(I+(L_1,L_2,L_3))} & : & 1 \ 1 \ 1  & \Rightarrow & r_3(R/I) = 2 & \hbox{(by
WLP)}

\end{array}
\]
The first use of WLP came because we had a complete intersection; the second came
because we were in a ring with two variables.  So we see that having WLP is extremely
helpful in computing the reduction numbers. 
\qed
\end{example}


\section{Setting the stage} \label{set stage}

The main goal of this paper is to sketch a progression of results on
Hilbert functions for sets of points, starting in $\mathbb P^2$ (mostly due to
Davis, cf.\ 
\cite{davis}), then moving to some of the generalizations of these results to higher
projective space obtained by Bigatti, Geramita and myself (cf.\ \cite{BGM}), and
finally recent extensions of some of these results obtained by Ahn and myself
(cf.\ \cite{AM}). In Section \ref{general results} we will focus on the general case,
while in Section
\ref{UPP results} we will specialize to the case of UPP. For simplicity of exposition
we will usually focus on  reduced schemes, but many of  these results extend to the
non-reduced case.

Let $Z$ be a zero-dimensional scheme.  The central assumption in this paper
will be the following:

\begin{equation} \label{main assumption}
\Delta h_Z(d) = \Delta h_Z(d+1) = s, \hbox{ for some $d, s$.}
\end{equation}  In both Section \ref{general results} and Section \ref{UPP results},
we will be focusing on two situations:

\begin{itemize}

\item[(a)]
$d \geq s$ (which was used in \cite{BGM}) 

\item[(b)] $d > r_2(R/I_Z)$ (which was used in
\cite{AM}).  
\end{itemize}

We will see that (a) corresponds to a certain kind of maximal growth (see Remark
\ref{s = s max gr}), while (b) in general does not (see Remark \ref{connection betw
hyp}).  This is the striking aspect of the results centered around (b) (see
especially Theorem \ref{AM general} and Theorem
\ref{AM UPP results}), that strong results can be obtained even without the power of
the Gotzmann machinery behind us. In this section we make a series of remarks in
preparation for the discussion in the coming sections.

\begin{remark} \label{comment on 334}
In Example \ref{334 example} we have (\ref{main assumption}), but neither $d
\geq s$ nor $d > r_2(R/I_Z)$ holds!  (Neither do the conclusions that we will
mention in the coming sections.)
\qed
\end{remark}

\begin{remark} \label{connection betw hyp}

What is the connection between (a) and (b) above?  If $d \geq s$ then it can be shown
that $d > r_2(R/I_Z)$ (see \cite{AM} Remark 4.4).  But in general
$d > r_2(R/I_Z)$ does  {\em not} imply $d \geq s$, and indeed it can happen that 
$d > r_2(R/I_Z)$ but $d$ is much smaller than $s$. 
In Example \ref{334 example}, $s=8$ and $r_2(R/I_Z) = 3$.  \qed
\end{remark}

\begin{remark} \label{max gr in p2} 
We now remark that in $\mathbb P^2$, condition (\ref{main assumption}) essentially
always corresponds to case (a) (there is only one, very special, exception).  Assume
that $Z \subset \mathbb P^{2}$ is a zero-dimensional scheme.  For $d \leq \alpha-1$
we have $\Delta h_Z(d) = d+1$ and for $d \geq \alpha$ we have that $\Delta h_Z$ is
non-increasing, so in particular  $\Delta h_Z(d) \leq \alpha$.

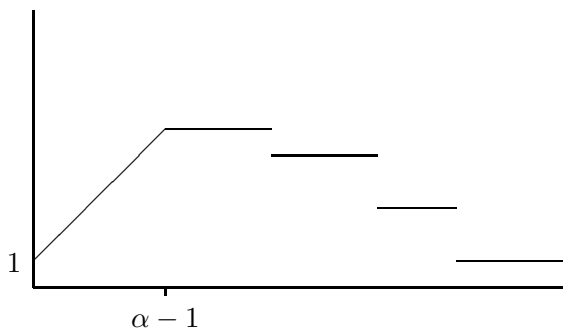
\vbox{$$
\begin{picture}(130,140)(-20,-100)
\put(-30,15){\line(0,-1){105}}
\put(-30,-90){\line(1,0){200}}
\put(-30,-80){\line(1,1){50}}
\put(-40,-84){\small $1$}
\put(20,-30){\line(1,0){40}}
\put(7,-105){\small $\alpha-1$}
\put(20,-90){\line(0,-1){3}}
\put(60,-40){\line(1,0){40}}
\put(100,-60){\line(1,0){30}}
\put(130,-80){\line(1,0){40}}
\end{picture}$$}

\bigskip

Of course there may not be any flat part at all, but the point is that the function
cannot increase past degree $\alpha-1$.
In any case we have
\[
[\Delta h_Z(d)  = \Delta h_Z(d+1) = s] \Rightarrow 
\left \{
\begin{array}{ll}
d = s-1, & \hbox{if $d = \alpha-1$ and $s=\alpha$} \\
d \geq s, & \hbox{if $d \geq \alpha$}
\end{array}
\right.
\]
In the first case in the above calculation, we  have that
$\dim(I_Z)_\alpha = 1$.  Notice that in this case $I_Z$ is zero in degree $d =
\alpha-1$, but $(I_Z)_{d+1} = (I_Z)_\alpha$ defines a hypersurface (or equivalently a
curve) in
$\mathbb P^2$.  \qed
\end{remark}

\begin{remark} \label{s = s max gr}
 Assume that $Z \subset \mathbb P^{n}$ is a zero-dimensional scheme.  
If $d \geq s$ then the condition $\Delta h_Z (d) = \Delta h_Z(d+1) = s$ implies that
the growth of the Artinian reduction from degree $d$ to degree $d+1$ is maximal in
the sense of Definition \ref{o-seq}.  Indeed, the $d$-binomial expansion of $s$ is
\[
s = \binom{d}{d} + \dots + \binom{d-s+1}{d-s+1}
\]
so
\[
s^{\langle d \rangle} = \binom{d+1}{d+1} + \dots + \binom{d-s+2}{d-s+2}=s
\]
as claimed.  \qed
\end{remark}

\begin{remark} \label{delta max not imply I max}
In practice, it usually is the case that $\Delta h_Z$ having maximal growth from
degree $d$ to degree $d+1$ does {\em not} imply that $h_Z$ itself has maximal growth
from degree $d$ to degree $d+1$.  If it did, Gotzmann's results (Theorem
\ref{gotzmann persistence} and Remark 
\ref{interpret max gr}) would immediately apply to $h_Z$.  The interest in the
results described below is that similar powerful conclusions come from this maximal
growth of the first difference (see Proposition \ref{BGM saturation result} and
Remark \ref{comparison}).  And as we will see, even more surprising is that similar
results can be deduced at times even when the first difference does not have maximal
growth. 
\qed
\end{remark}


\section{General results} \label{general results}

A by-now classical result of Davis (cf.\ \cite{davis}) is the following.  Note that
there is no uniformity assumption on $Z$.  In the next section we will discuss the
refinements that are possible when we assume UPP.

\begin{thm}[Davis] \label{davis thm}
Let $Z \subset \mathbb P^{2}$ be a  zero-dimensional scheme.  If $\Delta h_Z(d) =
\Delta h_Z(d+1) = s$ for $d \geq s$, then $(I_Z)_d$ and $(I_Z)_{d+1}$ both have a
GCD,
$F$, of degree $s$.  If $Z$ is reduced then so is $F$.  
The polynomial $F$ defines a $\left \{ 
\begin{array}{l}
\hbox{hypersurface} \\
\hbox{curve}
\end{array} \right. $
in $\mathbb P^{2}$.  If $Z_1 \subset Z$ is the subscheme of $Z$ lying on $F$ (defined
by
$[I_Z + (F)]^{sat}$) and $Z_2$ is the ``residual'' scheme defined by $[I_Z :F]$,
then there are formulas relating the Hilbert functions of $I_Z, I_{Z_1}$ and
$I_{Z_2}$.
\end{thm}

\begin{remark} \label{iarrobino comment}
It is worth mentioning that the Artinian version of Theorem \ref{davis thm} was
proved earlier, by A.\ Iarrobino (cf.\ \cite{iarrobino artin}, page 56).  \qed
\end{remark}

The paper \cite{BGM} began with the
observation that Davis' result, Theorem \ref{davis thm}, was really about maximal
growth of the function $\Delta h_Z(t)$ from degree $d$ to degree $d+1$, thanks to
Remark \ref{max gr in p2} and Remark \ref{s = s max gr} above.  These remarks show
that in $\mathbb P^2$, $\Delta h_Z$ can have two different kinds of maximal growth:
either $h_Z$ takes on the value of the polynomial ring (before the initial degree of
the ideal), or else $\Delta h_Z$ takes the constant value $s \geq d$.  This latter
is the only interesting case.  

Notice that since $Z$ is a zero-dimensional scheme, eventually the Hilbert
function is constant, so eventually $\Delta h_Z$ is zero.  So the condition that
$\Delta h_Z(d) = \Delta h_Z(d+1) = s$ directly gives us information not so much
about $I_Z$ as about $(I_Z)_{\leq d}$, the ideal generated by the components of
degree $\leq d$. We have to deduce properties of $I_Z$ from this, via Proposition
\ref{BGM saturation result}.  Gotzmann's results tell us that in a general Artinian
reduction of $R/I_Z$ (which has Hilbert function $\Delta h_Z$), the degree $d$
component of the ideal  $J = \frac{(I_Z,L)}{(L)}$ agrees with the degree $d$
component of the saturated ideal of a subscheme of degree $s$ in the line defined by
$L$ (cf.\ Remark \ref{interpret max gr}).  Since this corresponds to a hyperplane
section by $L$, that means that the base locus of the linear system $|I_d|$ contains
a curve of degree $s$, which is the GCD in question.  The results about the
subscheme of $Z$ lying on this GCD and the ``residual" subscheme come from some
algebraic arguments chasing exact sequences.

As noted in the way we phrased Davis'
theorem, this GCD can be viewed as a curve, or as a hypersurface in $\mathbb P^2$. 
These correspond, respectively, to the interpretations (that happen to coincide in
this case) that the maximal growth for $\Delta h_Z$  takes the constant value $s \leq
d$, and it takes the largest possible growth short of being equal to the Hilbert
function of the whole polynomial ring itself.

In \cite{BGM} we extended this to $\mathbb P^{n-1}$, and we called these kinds of
maximal growth {\em growth like a curve},  and {\em growth like a hypersurface}, 
respectively.  We proved many results for $Z \subset \mathbb P^{n-1}$.  In
particular, we showed that viewing $F$ in Theorem \ref{davis thm} as a curve, and
viewing it as a  hypersurface,  both extend in separate directions when we
move to higher projective spaces.  

In this paper we focus on the former.  We now give a summary of the results in the
``curve" direction found in \cite{BGM} that do not assume UPP.

\begin{thm}[\cite{BGM}] \label{BGM general results}
Let $Z \subset \mathbb P^{n-1}$ be a {reduced} zero-dimensional scheme.  Assume that
$\Delta h_Z(d) = \Delta h_Z(d+1) = s$ for some $d \geq s$.   Then

\begin{itemize}
\item[(a)] $\langle (I_Z)_{\leq d} \rangle$ is the {\bf saturated} ideal of
a curve $V$ of degree $s$.  $V$ is not necessarily unmixed, but it is
{\bf {reduced}} and {\bf $d$-regular}. 
\end{itemize}

Let $C$ be the unmixed one-dimensional part of $V$.  Let $Z_1$ be the subscheme of
$Z$ lying on $C$ (defined by $[I_Z + I_C]^{sat}$) and $Z_2$ the residual scheme
(defined by $[I_Z : I_C]$).  Then

\begin{itemize}
\item[(b)] $\langle (I_{Z_1})_{\leq d} \rangle = I_C$ and $I_C$ is
\underline{$d$-regular}.  (Part (a) was for $V$.  This one is not surprising, since
$V$ consists of $C$ plus some points.)

\item[(c)] There are formulas relating the Hilbert functions.

\end{itemize}

\end{thm}

How might these results be improved?  Recall that $\Delta h_Z(d) =$ \linebreak 
$\Delta h_Z(d+1) = s$ for $d \geq s$ guarantees that $\Delta h_Z$ has maximal growth
from degree $d$ to degree $d+1$.  

\begin{enumerate}

\item Weaken the condition $\Delta h_Z (d) = \Delta h_Z(d+1) = s$, e.g.\ to the
condition $\Delta h_Z(d) = \Delta h_Z(d+1) +1$ (possibly even maintaining the
assumption $d \geq s$).  It may be that at least ``usually" something similar will
hold.  But there will be important differences.  This is still an open direction,
although in her thesis Susan Cooper (cf.\ \cite{cooper}), a student of Tony
Geramita,  is working on questions related to this.

\item Weaken the assumption $d \geq s$.  This is the approach we take.  We will
assume only $d > r_2(R/I_Z)$.  This is a beautiful idea, conceived by my co-author,
Jeaman Ahn.  It is very striking to see how much carries over to this case.
\end{enumerate}

\begin{thm}[\cite{AM}] \label{AM general}
Let $Z \subset \mathbb P^{n-1}$ be a reduced zero- \linebreak dimensional scheme. 
Assume that
$\Delta h_Z(d) = \Delta h_Z(d+1) = s$ for some $d > r_2(R/I_Z)$.   Then

\begin{itemize}
\item[(a)] $\langle (I_Z)_{\leq d} \rangle$ is the {\underline{saturated}} ideal of
a curve $V$ of degree $s$.  $V$ is not necessarily unmixed  or reduced, but it
is $d$-regular. 
\end{itemize}

Let $C$ be the unmixed one-dimensional part of $V$.  Let $Z_1$ be the subscheme of
$Z$ lying on $C$ (defined by $[I_Z + I_C]^{sat}$) and $Z_2$ the residual scheme
(defined by $[I_Z : I_C]$).  Then

\begin{itemize}
\item[(b)] $\langle (I_{Z_1})_{\leq d} \rangle = I_C$ and $C$ is $d$-regular.

\item[(c)] There are formulas relating the Hilbert functions.

\item[(d)] {\bf If we also assume that} $h^1({\mathcal I}_{C_{red}}(d-1)) = 0$ then
$V$ is reduced and $C = C_{red}$ is $d$-regular.

\end{itemize}

\end{thm}

\begin{remark}
 We stress that in (a), we do not necessarily obtain that $V$ is reduced.  This is
the first new twist.  If $V$ is not reduced, then the top dimensional part, $C$, is
not necessarily reduced.  This curve $C$ is $d$-regular.  If it is not reduced,
however, it is supported on a reduced curve, $C_{red}$.  Surprisingly, $C_{red}$
being $d$-regular does not follow from $C$ being $d$-regular, and there are
examples where it is not true.  But with the extra assumption, (d) delivers this
conclusion.  Example \ref{not decr type} below illustrates what can happen.  \qed
\end{remark}

Theorem \ref{AM general} combines different results of \cite{AM}.  The proofs heavily
use generic initial ideals and results of Green, Bayer, Stillman, Galligo, etc.

We end this section with a result that incorporates WLP and also higher differences
of the Hilbert function and higher reduction number.

\begin{thm}[\cite{AM} Theorem 6.7] \label{delta 2}
Let $Z$ be a zero-dimen\-sional subscheme of $\mathbb P^{n-1}$,
$n>3$, with \textup{WLP}. Suppose that
\[
\Delta^2 h_Z (d)=\Delta^2 h_Z (d+1)=s
\]
for $r_2(R/I_Z)>d>r_3(R/I_{Z})$. Then
$\langle (I_Z)_{\leq d}\rangle$ is a saturated ideal defining a
two-dimensional subscheme of degree $s$ in $\mathbb P^{n-1}$, and it
is $d$-regular.
\end{thm}


\section{Results on Uniform Position} \label{UPP results}

We now investigate the effects of assuming that our points have UPP.  We begin again
in $\mathbb P^2$.  

\begin{thm} \label{P2 UPP results}
Let $Z \subset \mathbb P^2$ be a reduced set of points with UPP.  Then
we have

\begin{itemize}
\item[(a)] {\rm (\cite{ger-mar}, \cite{MR})} The component of $I_Z$ of least degree,
$\alpha$, contains an irreducible form (hence the general such form is irreducible). 
In particular, this holds if there is only one such form, up to scalar multiple. 
This is in fact true not only in $\mathbb P^2$ but also in $\mathbb P^{n-1}$.

\item[(b)] {\rm (\cite{harris})} $\Delta h_Z$ is of {\underline{decreasing
type}}, i.e.
\[
\hbox{if } \Delta h_Z(d) > \Delta h_Z(d+1) \hbox{ then } 
\Delta h_Z(t) > \Delta h_Z(t+1)
\]
 for all $t \geq d$ as long as $\Delta h_Z(t) > 0$.

\item[(c)]  If $\Delta h_Z(d) = \Delta h_Z(d+1) = s$ then $s =
\alpha$ (see Remark \ref{max gr in p2}), $(I_Z)_t = (F)_t$ for all $t \leq d+1$, and
the points of $Z$ all lie on the irreducible curve defined by $F$.

\end{itemize}
\end{thm}

As before, this theorem was also extended in [Bigatti-Geramita-Migliore] to
$\mathbb P^{n-1}$ for the case of ``maximal growth like a hypersurface," but here we
focus on extending it for ``maximal growth like a curve."  (See the discussion
preceding Theorem \ref{BGM general results}.)  In the following we repeat part of
Theorem
\ref{BGM general results} and underline the new features.

\begin{thm}[\cite{BGM}] \label{BGM UPP results}
Let $Z \subset \mathbb P^{n-1}$ be a {reduced} zero-dimensional scheme with UPP. 
Assume that $\Delta h_Z(d) = \Delta h_Z(d+1) = s$ for some $d \geq s$.   Then

\begin{itemize}
\item[(a)] $\langle (I_Z)_{\leq d} \rangle$ is the saturated ideal of
a curve $V$ of degree $s$.  $V$ is reduced, $d$-regular, \underline{unmixed} and
\underline{irreducible}.

\item[(b)] Since $V$ is unmixed, its top dimensional part $C$ is equal to
$V$.  Hence $Z \subset C$.

\item[(c)] $\langle (I_{Z})_{\leq d} \rangle = I_C$ and $I_C$ is
{$d$-regular}. 
\end{itemize}

\end{thm}

\begin{remark}
The hardest part is irreducibility, and it strongly uses the assumption $d \geq s$. 
\qed
\end{remark}

Still under the hypothesis of \cite{BGM} that $d \geq s$, the following extension of
``decreasing type" was observed in [Ahn-Migliore].  This completes the
picture of extending the $\mathbb P^2$ result to higher projective space under the
assumption of ``maximal growth like a curve."

\begin{corollary}[\cite{AM}] \label{decr type in Pn}
If $Z \subset \mathbb P^{n-1}$ has UPP and $\Delta h_Z(d) = $ \linebreak $\Delta
h_Z(d+1) >
\Delta h_Z(d+2)$ for some $d \geq s$, then $\Delta h_Z(t) > \Delta h_Z(t+1)$ for all
$t
\geq d+1$ as long as $\Delta h_Z(t) > 0$.
\end{corollary}

We now turn to the situation of \cite{AM}, where we only assume that $d >
r_2(R/I_Z)$.  It is somewhat surprising to see what we retain and what we lose. 
Compare this with Theorem \ref{AM general} (part of which is repeated here).

\begin{thm}[\cite{AM}] \label{AM UPP results}
Let $Z \subset \mathbb P^{n-1}$ be a reduced zero- \linebreak dimensional scheme with
UPP.  Assume that $\Delta h_Z(d) = \Delta h_Z(d+1) = s$ for some $d > r_2(R/I_Z)$.  
Then

\begin{itemize}
\item[(a)] $\langle (I_Z)_{\leq d} \rangle$ is the saturated ideal of an
\underline{unmixed} curve, $V$, of degree $s$, and it is $d$-regular.  $V$ is not
necessarily  reduced or irreducible.

\item[(b)] Since $V$ is unmixed, its top dimensional part $C$ is equal to
$V$.  Hence $Z \subset C$.

\item[(c)] $\langle (I_{Z})_{\leq d} \rangle = I_C$ and $I_C$ is
{$d$-regular}. 

\item[(d)] {\bf If we also assume that} $h^1({\mathcal I}_{C_{red}}(d-1)) = 0$ then
$V$ is reduced and $C = C_{red}$ is $d$-regular.

\end{itemize}

\end{thm}

\begin{remark}
Why might we have hoped that ``reduced and irredu\-cible'' would still hold in part
(a) of Theorem
\ref{AM UPP results}? Recall the following.

\begin{itemize}

\item We saw in Theorem \ref{P2 UPP results} that if $Z \subset \mathbb P^{n-1}$ is a
set of points in $\mathbb P^{n-1}$ with UPP, it is known that a general element of
smallest degree $\alpha$ is reduced and irreducible.  When this element is unique (up
to scalar multiple), say $F$, we have $\langle (I_Z)_{\leq \alpha}
\rangle$ is the saturated ideal of a reduced, irreducible hypersurface which is
$\alpha$-regular.  This is the base locus in degree $\alpha$, and just having UPP is
enough to guarantee the irreducibility. 

\item Recall  from Theorem \ref{BGM UPP results} that if $d \geq s$ then it {\em is}
true that the base locus, $C$, is reduced and irreducible.

\end{itemize}

In the next section we will give  examples to show that in fact ``reduced and
irreducible" does not necessarily hold (as opposed to merely being a gap in the
theorem).  \qed
\end{remark}


\section{Examples} \label{example section}

These examples will serve to illustrate that some of our results above are apparently
close to optimal.  We omit some details, and refer the reader to \cite{AM}.  We will
see that

\begin{enumerate}

\item In contrast to Corollary~\ref{decr type in Pn}, for a set of points in $\mathbb
P^{n-1}$ with UPP satisfying $\Delta h_Z(d) = \Delta h_Z(d+1) > \Delta h_Z(d+2)$ for
$d > r_2(R/I_Z)$ (instead of $d \geq s$), it is not necessarily true that $\Delta
h_Z$ is strictly decreasing beyond this point.  This is shown in Example~\ref{not
decr type}.

\item In  Theorem \ref{AM UPP results}, when we assume only $d > r_2(R/I_Z)$, it is
in fact  true (as claimed) that $V$ is not necessarily reduced or irreducible.  This
stands in contrast to Theorem \ref{BGM UPP results}, when we assumed $d \geq s$. 
This is  shown in Example \ref{not decr type}.

\item Upon realizing that in Theorem \ref{AM UPP results} $V$ is not necessarily
reduced and irreducible, one might hope that at least the points of $Z$ are
restricted to one component of $V$ that {\em is} reduced and irreducible.  Examples
\ref{not irred} and \ref{two sextics} show that even this is not true.

\item Examples \ref{can't extend} and \ref{back to ex} illustrate some of the
difficulties of extending these results past the case $d > r_2(R/I_Z)$, showing that
these results are optimal in a sense.

\end{enumerate}

We begin with an example of Chardin and D'Cruz.  It
settles an old question of independent interest.

\begin{example}[\cite{CD}] \label{CD example}
Consider the family of complete intersection \linebreak ideals 
\[
I_{m,n} := (x^mt - y^mz,
z^{n+2}-xt^{n+1}) \subset k[x,y,z,t].
\]
  Then for $m,n \geq 1$, $\reg(I_{m,n}) =
m+n+2$ while $\reg(\sqrt{I_{m,n}}) =$ \linebreak $  mn+2$.  Hence the regularity of
the radical may be much larger than the regularity if the ideal itself!!  In
particular, if we take $m=n=4$ then we obtain the following:  $I_{4,4}$ is a
complete intersection of type $(5,6)$.  As such, it has degree 30, and has a Hilbert
function whose first difference is

\begin{center}
\begin{tabular}{l|cccccccccccccccccccccccccccccc}
deg & 0 & 1 & 2 & 3 & 4 & 5 & 6 & 7 & 8 & 9 & 10 & 11 & 12 \\
$\Delta h$ & 1 & 3 & 6 & 10 & 15 & 20 & 24 & 27 & 29 & 30 & 30 & \dots
\end{tabular}

\end{center}

On the other hand, $\sqrt{I_{4,4}}$ can be computed on a computer
program, e.g.\
{\tt macaulay} (cf.\ \cite{macaulay}).  It has degree 26 (hence $I_{4,4}$ is not
reduced), and its Hilbert function has first difference

\medskip

{\small

\begin{tabular}{l|cccccccccccccccccccccccccccccc} 
deg & 0 & 1 & 2 & 3 & 4 & 5 & 6 & 7 & 8 & 9 & 10 & 11 & 12 & 13 & 14 & 15 \\ \hline 
$\Delta h$ & 1 & 3 & 6 & 10 & 15 & 20 & 24 & 27 & 29 & 29 & 29 & 29 & 28 & 28 &
28 & 27 \end{tabular}

\medskip

\begin{tabular}{l|cccccccccccccccccccccccccccccc} 
deg & 16 & 17 & 18 & 19 & 20 & 21 \\ \hline
$\Delta h$ & 27 & 27 & 26 & 26 & 26 & 26 & \dots  
\end{tabular}
}

\medskip

In \cite{CD}, the authors were not so concerned with the geometry of this curve.  We
will modify this slightly, taking a more geometric approach, in Example \ref{not
decr type}. 
\qed
\end{example}

\begin{remark}
Example \ref{CD example} shows that a curve $C \subset \mathbb P^3$ can have the
property that $\reg(\sqrt{I_C})$ may be (much) larger than $\reg(I_C)$.  However, it
seems to still be an open question whether this can happen, for instance, for a
smooth curve.  Chardin and D'Cruz also study the surface case, where other
interesting phenomena occur.  For curves, see also Ravi \cite{ravi}.  \qed
\end{remark}

\begin{example}[\cite{AM} Example 5.7] \label{not decr type}
In this example we recall that $k$ has characteristic zero.  
We will be basing our example on the case $m = n = 4$ of the example of Chardin and
D'Cruz.  It is obtained by taking a geometric interpretation.  Let $R =
k[x,y,z,t]$.  

 Let $I_\lambda = (z,t)$. Let $F \in (I_\lambda)_5$ be a homogeneous polynomial that
is  smooth along $\lambda$.  Let $I' = I_\lambda^5 + (F)$.  $I'$ is the saturated
ideal of a non-ACM curve of degree 5 corresponding to the divisor $D := 5 \lambda$
on $F$. In particular, viewed as a subscheme of $\mathbb P^3$, $D$ has degree 5.

Now, choose a general element $C$ in the linear system $|6H-D|$ on $F$.  $C$ is
smooth and irreducible.  Furthermore, $C$ has degree 25 (by liaison) and Hilbert
function with first difference 

{\small

\begin{tabular}{l|cccccccccccccccccccccccccccccc} 
deg & 0 & 1 & 2 & 3 & 4 & 5 & 6 & 7 & 8 & 9 & 10 & 11 & 12 & 13 & 14 & 15  \\
\hline 
$\Delta h$ & 1 & 3 & 6 & 10 & 15 & 20 & 24 & 27 & 29 & 29 & 29 & 29 & 28 & 28 &
28 & 27 \end{tabular}

\medskip

\begin{tabular}{l|cccccccccccccccccccccccccccccc} 
deg & 16 & 17 & 18 & 19 & 20 & 21 \\ \hline
$\Delta h$ & 27 & 27 & 26 & 26 & 26 & 25 & \dots  
\end{tabular}
}

We now let $Z$ consist of sufficiently many points of $C$, chosen generally, so
that $I_Z$ agrees with $I_C$ up to and including degree 21. 
$Z$ has UPP (since $C$ is smooth), and one checks that $r_2(R/I_Z) = r_2(R/I_C) = 
8$.  Taking
$d=9$, we see that Theorem \ref{AM UPP results} applies.  We obtain that $\langle
(I_Z)_{\leq 9} \rangle$ is the saturated ideal of an unmixed  curve
$V$  of degree 29 consisting of the union of $C$ and a subcurve of
$D$ of degree 4 supported on $\lambda$, hence $V$ is neither irreducible nor reduced.
Taking $d = 12$, $d = 15$, $d = 18$ slices away at the non-reduced part, and taking
$d = 21$ gives just $C$. 

A computation on {\tt macaulay} reveals that in fact the regularity of $I_C$ is
21. On the other hand, the curve $V$ of degree 29 has ideal $I_V$ with regularity  9.
Actually, it is worth recording the Betti diagram of $I_C$ and $I_V$, because they
display a surprising pattern:  For $I_V$ we have

\medskip

{\scriptsize 
\begin{verbatim}
                        total:      1     3     2 
                        --------------------------
                            0:      1     -     - 
                            1:      -     -     - 
                            2:      -     -     - 
                            3:      -     -     - 
                            4:      -     1     - 
                            5:      -     1     - 
                            6:      -     -     - 
                            7:      -     -     - 
                            8:      -     1     2 
\end{verbatim}
}
\medskip

\noindent while for $I_C$ we have

\medskip

{\scriptsize 
\begin{verbatim}
                     total:      1     7    10     4 
                     --------------------------------
                         0:      1     -     -     - 
                         1:      -     -     -     - 
                         2:      -     -     -     - 
                         3:      -     -     -     - 
                         4:      -     1     -     - 
                         5:      -     1     -     - 
                         6:      -     -     -     - 
                         7:      -     -     -     - 
                         8:      -     1     2     - 
                         9:      -     -     -     - 
                        10:      -     -     -     - 
                        11:      -     1     2     1 
                        12:      -     -     -     - 
                        13:      -     -     -     - 
                        14:      -     1     2     1 
                        15:      -     -     -     - 
                        16:      -     -     -     - 
                        17:      -     1     2     1 
                        18:      -     -     -     - 
                        19:      -     -     -     - 
                        20:      -     1     2     1 
\end{verbatim}
}
\medskip

\noindent Notice that the diagram for $I_V$ is a subdiagram of the diagram
for $I_C$, and that there is a striking simplicity to the diagram for $I_C$.  Notice
also that $V$ is ACM, while $C$ is not.
\qed

\end{example}

 Note that in the previous example, the points of $Z$  all lie on one reduced,
irreducible curve ($C$).  It is simply the case that $V$ contains another component,
which happens to also be non-reduced.  However, now we give another example to show
that it is not necessarily true that all the points lie on one irreducible component
of the base locus.

\begin{example}[\cite{AM} Example 5.8] \label{not irred}
In Example \ref{not decr type}, instead of choosing ``sufficiently many"
points on $C$,
instead choose $Z$ to consist of 192 general points on $C$ and one
general point of
$\lambda$.  The first difference of the Hilbert function of $Z$ is

\begin{center}
\begin{tabular}{l|cccccccccccccccccccccccccccccc}
deg & 0 & 1 & 2 & 3 & 4 & 5 & 6 & 7 & 8 & 9 & 10 & 11 &   \\
\hline $\Delta h_Z$ & 1 & 3 & 6 & 10 & 15 & 20 & 24 & 27 & 29 & 29 &
29 & 0
\end{tabular}

\end{center}

Note that this is exactly the same as what we would have if we had taken 193
general points of $C$.  Again, $r_2(R/I_Z) = 8$.  The base locus of
$(I_C)_{9}$ and $(I_C)_{10}$ is exactly the non-reduced and reducible curve of degree
29 mentioned in Example~\ref{not decr type}, so the Hilbert function of $Z$ is the
truncation of the Hilbert function given above.  And by the general choice of the
points, this will continue to be true regardless of which subsets we take.  Hence $Z$
has UPP and satisfies $\Delta h_Z(d) = \Delta h_Z(d+1) = s$ for some $d >
r_2(R/I_Z)$, but not all of the points of $Z$ lie on one reduced and irreducible
component of the curve of degree $s$ obtained by our result (since one point lies on
the non-reduced component).  \qed
\end{example}

\begin{example} {\rm [Ahn-Migliore] Example 5.9]} \label{two sextics}
We produced  Example \ref{not irred} by taking almost all of the points on
$C$, and just one off of $C$ (in fact, it was on $\lambda$).  In this example we
show a surprising instance where half of the points of $V$ are on one irreducible
component and the other half are on another irreducible component, but we still have
UPP.  We again omit some of the technical details.

Let $Q$ be a smooth quadric surface in $\mathbb P^3$, and as usual by abuse of
notation we
use the same letter $Q$ to denote the quadratic form defining this surface.
Let $C_1$ be a
{\em general} curve on $Q$ of type $(1,15)$, and let $C_2$ be a {\em general}
curve on $Q$ of type $(15,1)$.  Hence both $C_1$ and $C_2$ are smooth rational
curves of degree 16, and $C := C_1 \cup C_2$ is the complete
intersection of $Q$
and a form of degree 16.  Note that $C$ is arithmetically Cohen-Macaulay, but
$C_1$ and $C_2$ are not.

It is not difficult to compute the Hilbert functions of these curves.
We record
their first differences (of course there is no difference between behavior of
$C_1$ and behavior of $C_2$; this is important in the argument given in \cite{AM}):

\bigskip

{\small

\begin{tabular}{l|ccccccccccccccccccccccccccccccccccccccccccc}
degree & 0 & 1 & 2 & 3 & 4 & 5 & 6 & 7 & 8 & 9 & 10 & 11 & 12 & 13 \\ \hline
$\Delta h_C$ & 1 & 3 & 5 & 7 & 9 & 11 & 13 & 15 & 17 & 19 & 21 & 23 & 25 &
27 \\
$\Delta h_{C_i}$ & 1 & 3 & 5 & 7 & 9 & 11 & 13 & 15 & 17 & 19 & 21 & 23 &
25 & 27 
\end{tabular}

\medskip

\begin{tabular}{l|ccccccccccccccccccccccccccccccccccccccccccc}
degree & 14 & 15 & 16 & 17 & 18 \\ \hline
$\Delta h_C$ & 29 & 31 & 32 & 32 & 32 \\
$\Delta h_{C_i}$ & 29 & 16 & 16 & 16 & 16 
\end{tabular} }

\medskip

\noindent We now observe:
\begin{enumerate}

\item These first differences (hence
the ideals themselves) agree through degree 14, and in fact the only generator
before degree 15 is $Q$.

\item \label{hf values} By adding these values, we see that $h_C(18) = 352$
and \linebreak
$h_C(15) = 241$.

\item Since $C$ and $C_i$ are curves, these values represent the Hilbert
functions of
$I_C + (L)$ and $I_{C_i} +(L)$ for a general linear form $L$.

\item \label{value of r2} $r_2(R/I_C) = 16$ since $C$ is an arithmetically
Cohen-Macaulay curve.
\end{enumerate}

\bigskip

Let $Z_1$ (respectively $Z_2$) be a general set of 176
points on $C_1$
(respectively $C_2$).  So $Z := Z_1 \cup Z_2$ is a set of 352
points whose Hilbert function agrees with that of $C$ through degree 18.  In
particular, we have $\Delta h_Z(17) = \Delta
h_Z(18) = 32 = \deg
C$.  Furthermore, $r_2(R/I_Z) = 16$ by our observation \ref{value of r2}.\ above. 
Hence Theorem \ref{AM UPP results} applies, and we indeed have  that the component of
$I_Z$ in degree 17 defines $C$.  However, $C$ is not irreducible.  In \cite{AM}
Example 5.9, it is shown that
$Z$ has the Uniform Position Property.  Thus there is no chance of showing that
all the points must lie on a unique irreducible component in Theorem
\ref{AM UPP results}, under our hypothesis that $d > r_2(R/I_Z)$ (as was done in
\cite{BGM} when $d \geq s$).

To show UPP, it is enough to show that the union of any choice of $t_1$ points of
$Z_1$ (i.e.\ $t_1$ general points of $C_1$ for $t_1 \leq 176$) and $t_2$
points of $Z_2$ (i.e.\ $t_2$ general points of $C_2$ for $t_2
\leq 176$) has the truncated Hilbert function.  For example, if $t_1 = 150$ and $t_2
= 160$ then we have to show that $\Delta H(R/I_Z)$ has values
\[
1 \ \ 3 \ \ 5 \ \ 7 \ \ 9 \ \ 11 \ \ 13 \ \ 15 \ \ 17 \ \ 19 \ \ 21 \
\ 23 \ \ 25
\ \ 27 \ \ 29 \ \ 31 \ \ 32 \ \ 22 \  \ 0.
\]
Notice that we know that some subset has this Hilbert function, by \cite{GMR}.
We have to show that {\em all} subsets have this Hilbert function.  See \cite{AM},
Example 5.9, for the details.  \qed
\end{example}

\begin{example} \label{can't extend}
We saw that the condition $d \geq s$ of [Bigatti-Geramita-Migliore] was
improved in
\cite{AM} to
$d > r_2(R/I_Z)$.  There was some loss in the strength of the result, but a
surprising amount of it did go through. One might wonder if the condition $d >
r_2(R/I_Z)$ can be further improved.  But in fact, we saw already in Example
\ref{334 example} and Remark \ref{comment on 334} that this is not the case.  

Similarly, we have the following examples:
\[
\begin{array}{l}
\hbox{7  general points in $\mathbb P^3$ have $h$-vector 1 \ 3 \ 3} \\
\hbox{16 general points in $\mathbb P^3$ have  $h$-vector 1 \ 3 \ 6 \ 6} \\
\hbox{30 general points in $\mathbb P^3$ have  $h$-vector 1 \ 3 \ 6 \
10 \ 10} \\
\hbox{etc.}
\end{array}
\]
In each case, $r_2(R/I_Z)$ has the ``expected" value, say $d$.  We have $\Delta h_Z
(d) = \Delta h_Z(d+1) $ for $d = r_2(R/I_Z)$, but clearly the base locus of $(I_Z)_d$
is not one-dimensional.  See Remark \ref{comparison} as well.

Similar examples can easily be found in higher projective spaces.  A different sort
of example can be found in Example \ref{334 example}. \qed
\end{example}

\begin{example} \label{back to ex}
We return to Example \ref{WLP in families} to see how the theorems mentioned above
apply to that example.  Recall that we have two ideals, $I_1$ and $I_2$, both with
$h$-vector
$(1,3,6, 9, 11,11,11)$.  Taking $s=11$, it is clear that the results from \cite{BGM}
(Theorem \ref{BGM general results} and Theorem \ref{BGM UPP results}) do not apply
because of the hypothesis $d \geq s$.  As for the results of \cite{AM} (Theorem
\ref{AM general} and Theorem \ref{AM UPP results}), one can check that $r_2(R/I_1) =
4$ while $r_2(R/I_2) = 5$.  Hence both results from \cite{AM} apply to $I_1$, taking
$d = 5$, but not to $I_2$.  And indeed, we have seen that $\langle (I_1)_{\leq 5}
\rangle$ is the saturated ideal of the curve $C$ described in that example.  However,
it can easily be checked that
$\langle (I_2)_{\leq 5} \rangle$ is saturated, but it defines a curve of degree 10
rather than 11 (and it is not unmixed, although the unmixed part is ACM; in fact it
is a complete intersection of type $(2,5)$.). \qed
\end{example}

\end{document}